\documentclass[a4paper,12pt]{article}
\usepackage{tikz}
\usetikzlibrary{calc}
\usepackage{amssymb}
\usetikzlibrary{quotes,angles}
\usetikzlibrary{plotmarks}
\usetikzlibrary{arrows,shapes,positioning}
\usetikzlibrary{decorations.markings}
\usepackage{geometry}
\usepackage{caption}
\tikzstyle arrowstyle=[scale=2]
\tikzstyle directed=[postaction={decorate,decoration={markings,
    mark=at position .65 with {\arrow[arrowstyle]{stealth}}}}]
\tikzstyle reverse directed=[postaction={decorate,decoration={markings,
    mark=at position .65 with {\arrowreversed[arrowstyle]{stealth};}}}]
    \tikzstyle left directed=[postaction={decorate,decoration={markings,
    mark=at position -.62 with {\arrow[arrowstyle]{stealth}}}}]

\tikzstyle left reverse directed=[postaction={decorate,decoration={markings,
    mark=at position -.62 with {\arrowreversed[arrowstyle]{stealth};}}}]
\usepackage{indentfirst}
\usepackage[plainpages=false]{hyperref}
\usepackage{amsfonts,latexsym,rawfonts,amsmath,amssymb,amsthm,mathrsfs}
\usepackage{amsmath,amssymb,amsfonts,latexsym,lscape,rawfonts}

\usepackage[all]{xy}
\usepackage{eufrak}
\usepackage{makeidx}         
\usepackage{graphicx,psfrag}
\usepackage{array,tabularx}

\usepackage{setspace}

\newtheorem{thm}{Theorem}[section]
\newtheorem{cor}[thm]{Corollary}

\newtheorem{clm}[thm]{Claim}

\theoremstyle{remark}

\theoremstyle{definition}
\newtheorem{Def}[thm]{Definition}                                        %

\title{On partial uniqueness of complete non-compact Ricci flat metrics}

{\small{\author{Yuanqi Wang\thanks{University of Kansas, Lawrence, KS, USA.\ yqwang@ku.edu.}}}

\date{\vspace{-5ex}}

\begin{document}
\maketitle
\begin{abstract} Using  techniques for Caccioppoli inequality, on a fairly general class of complete non-compact K\"ahler manifolds with sub-quadratic volume growth, we show uniqueness of bounded $C^{1,1}$ solution to  Monge-Ampere equation. This does not a priori require any decay of the solution.  
\end{abstract}
\section{Introduction}
In an arbitrary K\"ahler class on a closed K\"ahler manifold  with vanishing $c_{1}$, uniqueness  of Ricci flat metric  can be proved via  integration by parts \cite{Calabi,Yau}.  
 On a non-compact complete K\"ahler manifold,  we show that under sub-quadratic volume growth and other mild conditions, even without decay,  there is still an integration by parts argument for uniqueness of bounded solution to  Monge-Ampere equation. The technique follows Caccioppoli inequality on regularity for  elliptic  equations. For example, see \cite{BM,DT, GM, KS}. 

A sub-harmonic function on a compact manifold must be a constant. In contrast and in general, the same fails on Euclidean domains even if we prescribe  constant boundary value. For example, the  function $|X|^{2}-1$ is sub-harmonic and vanishes on the boundary of the unit ball centered at the origin. 
Via this viewpoint, a complete non-compact Riemannian manifold is in between:  it is not compact, but  it is complete as a metric space and has no boundary. Our result shows, under the volume growth condition, a sub-solution to  Monge-Ampere equation with constant density $1$ is still a constant.

\begin{thm}\label{thm}  Let  $(M, \omega)$ be a complete non-compact K\"ahler manifold of complex dimension $n\geq 2$ and sub-quadratic volume growth. Let  $f$ be  a  continuous  function bounded from above on $M$.

\begin{enumerate}\item   Suppose $\omega$ has strictly sub-quadratic volume growth.  Then any bounded $C^{1,1}$  sub-solution, super-solution, or solution  $\phi$ to   \begin{equation}\label{MA0}(\omega+i\partial \overline{\partial}\phi)^{n}=\omega^{n}
\end{equation} 
  is a constant.  
\item Suppose   $e^{f}-1$ is integrable i.e. $\int_{M}|e^{f}-1|\omega^{n}<+\infty$. Then any bouned $C^{1,1}$  solution $\phi$  to \begin{equation}\label{MA}
(\omega+i\partial \overline{\partial}\phi)^{n}=e^{f}\omega^{n}
\end{equation}  has  bounded Dirichlet energy i.e. 
$\int_{M}|\nabla_{\omega}\phi|^{2}\omega^{n}\leq c_{DE}<+\infty,$
where \begin{equation}\label{equ cD}c_{DE}=1000n^{4n+4}|\phi|_{C^{0}(M)}[(1+K)^{2n} c_{vol} + |e^{f}-1|_{L^{1}(M,\omega^{n})}].\end{equation}
\item  Suppose additionally that $\omega$  satisfies weak Neumann Poincare inequality on annulus. Then any bounded $C^{1,1}$ solution $\phi$ to \eqref{MA0}   is a constant. 
\end{enumerate}

\end{thm}
The constant \eqref{equ cD} might not be optimal, but is effective. The terms involved are defined below.  
\begin{Def}\label{Def MA} 

A \textit{bounded $C^{1,1}$  solution, sub-solution, or super-solution   to \eqref{MA}}  is a real-valued function $\phi$  with the following properties.
\begin{itemize}
\item $\phi$ is twice continuously differentiable  under the holomorphic (smooth) manifold structure (not necessarily with norm bound).
\item $\phi$ satisfies \eqref{MA},   $(\omega+i\partial \overline{\partial}\phi)^{n}\geq e^{f}\omega^{n}$,  or $(\omega+i\partial \overline{\partial}\phi)^{n}\leq e^{f}\omega^{n}$ respectively.
\item There is a positive number $K\ (<\infty)$   such that
\begin{equation}\label{equ C11 condition} |\phi|_{C^{0}(M)}=\sup_{M}|\phi|<\infty
\ \textrm{and}\ 0< \omega_{\phi}\leq K \omega\  \textrm{on the whole}\ M. \end{equation}
\end{itemize}
Fix a point $o \in M$  as  distance origin and center of balls.  We say  a complete non-compact K\"ahler metric $\omega$ has \textit{sub-quadratic volume growth},
  if there is  a positive sequence $\rho_{i}\rightarrow +\infty$ such that \begin{equation}\limsup_{i\rightarrow \infty} \frac{Vol[B(2\rho_{i})\setminus B(\rho_{i})]}{\rho_{i}^{2}}< +\infty. 
\end{equation}

The value of the existing  limit superior is denoted by $c_{vol}$.
We say such an $\omega$ has \textit{strictly sub-quadratic volume growth} if $c_{vol}=0$ i.e.  \begin{equation}\label{equ vol slower than quadratic}
\lim_{i\rightarrow \infty}\frac{Vol[B(2\rho_{i})\setminus B(\rho_{i})]}{\rho_{i}^{2}}=0. 
\end{equation}

On the other hand, we say  it satisfies  \textit{weak Neumann  Poincare inequality on annulus} if there is $\underline{\rho}_{i}\rightarrow\infty$ and  $\mu_{i}$, 
 such that $\underline{\rho}_{i}\leq \rho_{i}$ when $i$ is large, and 
\begin{equation}\label{equ Poincare inequ} \int_{B(2\rho_{i})\setminus B(\rho_{i})}|\phi -\mu_{i}|^{2}\leq c_{P}\rho^{2}_{i}\int_{\{r\geq \underline{\rho}_{i}\}}|\nabla\phi|^{2},
\end{equation}
where $c_{P}$ is independent of $i$ or $\phi$, as long as $\phi$ is twice continuously differentiable.  Our argument \eqref{equ energy bd using Poincare} below is independent of $\mu_{i}$. 
\end{Def}
Suppose we have two solutions $\phi_{1}$ and $\phi_{2}$ to the general  volume form equation \eqref{MA}. Fo uniqueness, as long as the conditions hold,  we can apply Theorem \ref{thm}.1 or \ref{thm}.3   with reference metric being $\omega+i\partial \overline{\partial}\phi_{1}$ or $\omega+i\partial \overline{\partial}\phi_{2}$,  and  $\phi$ being $\pm(\phi_{2}-\phi_{1})$ respectively.

Theorem \ref{thm} partially addresses the uniqueness of  Tian-Yau   solutions \cite[Theorem 1.1]{TianYau}.  This particular  result  is  under  sub-quadratic volume growth,  and their solution  is bounded $C^{1,1}$.  
Let ``unique...(up to constant)" abbreviates  ``unique...up to addition by a real constant". 
\begin{cor}\label{Cor} Let $2>\alpha\geq 1$.  If the $(K,2,\beta)-$polynomial growth condition is strengthened to $(K,\alpha,\beta)$, the Tian-Yau solution $\varphi$ in \cite[Theorem 1.1]{TianYau} is the unique bounded $C^{1,1}$
solution (up to constant) to the Monge-Ampere equation  \cite[(1.1)]{TianYau}. The solution $u$ in  Hein's version \cite[Prop 4.1]{HeinThesis}, under $SOB(\beta)-$condition, $\beta\leq 2$,  is  the unique bounded $C^{1,1}$
solution (up to constant)  to the Monge-Ampere equation therein. Consequently, the following holds.

\begin{itemize}\item The Tian-Yau Ricci-flat space \cite[Theorem 4.1]{TianYau} of volume growth $O(r^{\frac{2n}{n+1}})$ 
 is the unique solution (up to constant) to the Monge-Ampere equation \cite[(1.1)]{TianYau} with reference metric $\omega_{N}$ \cite[(4.4)]{TianYau}. 
 
\item  
 All the gravitational instantons  in Hein's construction \cite[Theorem 1.5]{HeinJAMS}  are unique solutions (up to constant) to corresponding Monge-Ampere equations in  \cite[Prop 4.1]{HeinThesis}. 
\end{itemize}
\end{cor}

Our sub-quadratic volume growth is a sequential condition, and  is apparently implied by 
that ball of radius $R$ (centered at the base point) has volume  $\leq CR^{2}$, for large $R$ cf. \cite[Definition 1.1]{TianYau}.
Suppose the   weak Poincare inequality \eqref{equ Poincare inequ} on annulus  is implied by   $(K,2,\beta)-$polynomial growth and  other conditions  in \cite[Theorem 1.1]{TianYau}.  Then uniqueness of Tian-Yau solution  \cite[Theorem 1.1]{TianYau} holds in full generality. Nevertheless, for Hein's version  \cite[Prop 4.1]{HeinThesis}, we do have  Poincare inequality \cite[Prop 3.4]{HeinThesis}. Therefore uniqueness also holds for rigorous quadratic volume growth i.e. $SOB(2)-$case. 






Under faster than quadratic volume growth, uniqueness of Monge-Ampere solutions is implied by certain  decay on the K\"ahler  potential  $\phi$. See \cite[8.5 Theorem A4]{Joyce} for example. Geometric uniqueness of Ricci flat metrics, as in  \cite{CH,
CH2, HHN, HeinThesis, HeinJAMS, Joyce, Sze}, 
 usually involves both Monge-Ampere uniqueness and $i\partial\overline{\partial}-$lemma under decay conditions.  For recent work on Liouville theorem of Monge-Ampere equations on product manifolds, see  Hein \cite{HeinCPAM} for example. For recent work on Liouville theorem of metric Laplacians on certain non-compact complete manifolds, see Sun-Zhang  \cite{SZ} and Carron  \cite{Carron} for examples. For earlier work on rigidity of  harmonic functions on non-compact complete Riemannian manifolds, see Cheng-Yau \cite{CY} for example. 
\section{Proof}
\textbf{Convention}:  Unless otherwise specified, the metric for   gradient is  $\omega$. The integrals and  volumes are with respect to the top degree form $\omega^{n}$.
\subsection*{For Theorem \ref{thm}.1}
Write $\omega_{\phi}$ for $\omega+i\partial \overline{\partial}\phi$.
We have the difference
\begin{equation}\label{equ difference between two Kah forms}
\omega^{n}_{\phi}-\omega^{n}=i\partial \overline{\partial}\phi\wedge Q.
\end{equation}
where  $Q=\omega^{n-1}_{\phi}+\omega^{n-2}_{\phi}\wedge \omega+...+\omega^{n-1}$. By positivity and  $C^{1,1}-$condition \eqref{equ C11 condition},
we verify $$\omega^{n-1} \leq Q\leq n^{2n}(1+K)^{n}\omega^{n-1}.$$
Actually, any constant depending only on the data in Theorem \ref{thm} and Definition
\ref{Def MA} suffices, but we want explicit constant,   though  not necessarily optimal. 
 Moreover, 
\begin{equation}\label{equ positive difference}
 \left\{ \begin{array}{c}i\partial \overline{\partial}\phi\wedge Q\ \geq 0\ \textrm{for subsolution,  and} \\
i\partial \overline{\partial}\phi\wedge Q\ \leq 0\ \textrm{for super-solution}.
\end{array}\right. 
\end{equation}

Because $\phi$ is bounded,  
\begin{itemize}
\item if it is a sub-solution, by adding a constant if necessary, we assume $ \min\phi\geq 1$;

\item  if it is a super-solution, by adding a constant if necessary, we  assume $ \max\phi\leq -1$.
\end{itemize}
 Let $\chi$ be compactly supported Lipschitz function. We  multiply both hand sides in \eqref{equ positive difference}  by $\chi^{2}\phi$. In either case, because of the definite sign of $\chi^{2}\phi$, we find 
\begin{equation}\label{equ positive difference either case}
\chi^{2}\phi\cdot i\partial \overline{\partial}\phi\wedge Q \geq 0.
\end{equation}


We  integrate \eqref{equ positive difference either case} by parts: 
\begin{equation}\label{equ initial ibp}
\int_{M}\chi^{2}i\partial \phi\wedge \overline{\partial}\phi \wedge Q \leq -2\int_{M}\phi\chi \cdot i\partial \chi\wedge \overline{\partial}\phi\wedge Q.
\end{equation}
The left side of \eqref{equ initial ibp} is bounded from below by
\begin{equation*}
\int_{M}\chi^{2}i\partial \phi\wedge \overline{\partial}\phi \wedge Q \geq \int_{M}\chi^{2}i\partial \phi\wedge \overline{\partial}\phi \wedge \omega^{n-1} \geq \frac{1}{2n}\int_{M}\chi^{2}|\nabla \phi|^{2}.
\end{equation*}
On the other hand, Cauchy-Schwartz yields an upper bound on   the right side of \eqref{equ initial ibp}:
\begin{eqnarray*}
& &|-2\int_{M}\phi\chi i\partial \chi\wedge \overline{\partial}\phi\wedge Q|\leq 2n^{2n}(1+K)^{n}\int_{M}|\phi||\chi| |\nabla \chi| |\nabla \phi|
\\&\leq & 2n^{2n}(1+K)^{n} (\int_{M} \chi^{2} |\nabla \phi|^{2})^{\frac{1}{2}}(\int_{M} \phi^{2} |\nabla \chi|^{2})^{\frac{1}{2}}.
\end{eqnarray*}
The above two inequalities and \eqref{equ initial ibp} imply 
 \begin{equation}\label{equ energy bound with chi} \int_{M}\chi^{2}|\nabla \phi|^{2}\leq 4n^{2n+1}(1+K)^{n} (\int_{M} \chi^{2} |\nabla \phi |^{2})^{\frac{1}{2}}(\int_{M} \phi^{2} |\nabla \chi|^{2})^{\frac{1}{2}}.\end{equation}

Either  $\int_{M}\chi^{2}|\nabla \phi|^{2}$ is $0$ or not, we find
\begin{equation}\label{equ the preliminary cacciopolli}
 \int_{M}\chi^{2}|\nabla \phi|^{2}\leq 16n^{4n+2}(1+K)^{2n}\int_{M} \phi^{2} |\nabla \chi|^{2}.
\end{equation}

Now let $\chi$ be the following piece-wise linear  function in the distance $r$ from the base point $o$, and $\rho>0$.
\begin{equation}\chi(r)=\left\{ \begin{array}{c}1\ \ \textrm{when}\ \ r\leq \rho,\\
2-\frac{r}{\rho}\ \ \textrm{when}\ \ \rho\leq r\leq 2\rho,\\
0\ \ \textrm{when}\ r\geq 2\rho. 
\end{array}\right.
\end{equation}
Apparently, $\chi$ is Lipschitz. We find
\begin{eqnarray}
& & \int_{B(\rho)}|\nabla \phi|^{2}\label{equ stripes} \leq  \frac{16n^{4n+2}(1+K)^{2n}}{\rho^{2}}\int_{B(2\rho)\setminus B(\rho)} |\phi|^{2} \nonumber
\\&\leq & 16n^{4n+2}(1+K)^{2n}|\phi|_{C^{0}(M)} \cdot \frac{Vol[B(2\rho)\setminus B(\rho)]}{\rho^{2}}. 
\end{eqnarray}
Let $\rho=\rho_{i}$ and $i \rightarrow \infty$. Monotone convergence theorem implies $|\nabla \phi|^{2}$ is integrable on $M$ and 

\begin{equation} \label{equ energy <=0} \int_{M}|\nabla \phi|^{2}\leq 16n^{4n+2}(1+K)^{2n}|\phi|_{C^{0}(M)} \cdot \limsup_{i\rightarrow \infty} \frac{Vol[B(2\rho_{i})\setminus B(\rho_{i})]}{\rho_{i}^{2}}= 0. \end{equation}
This means $\phi$ is a constant. 

\subsection*{For Theorem \ref{thm}.2}
In this case we do not add any constant to $\phi$. It is a solution by assumption.  The equality still holds if we multiply $i\partial\overline{\partial}\phi\wedge Q=(e^{f}-1)\omega^{n}$ by $\chi^{2}\phi$. The same argument \eqref{equ difference between two Kah forms}---\eqref{equ energy bound with chi} with $e^{f}\omega^{n}$ instead of $\omega^{n}$ yields
\begin{eqnarray*}& & \int_{M}\chi^{2}|\nabla \phi|^{2}
\\&\leq & 4n^{2n+1}(1+K)^{n}(\int_{M} \chi^{2} |\nabla \phi |^{2})^{\frac{1}{2}}(\int_{M} \phi^{2} |\nabla \chi|^{2})^{\frac{1}{2}}+\int_{M}\chi^{2}|\phi (e^{f}-1)|
\\&\leq &4n^{2n+1}(1+K)^{n} (\int_{M} \chi^{2} |\nabla \phi |^{2})^{\frac{1}{2}}(\int_{M} \phi^{2} |\nabla \chi|^{2})^{\frac{1}{2}}+|\phi|_{C^{0}(M)}\cdot |e^{f}-1|_{L^{1}(M)}.
\end{eqnarray*}

The following elementary claim is proved by completing square. 
\begin{clm}Let $a,\ b,\ c,\ d$ be non-negative numbers such that 
$$a\leq 2d\sqrt{a}\cdot \sqrt{b}+c.$$
Then $a\leq 2d^{2}b+c+2\sqrt{d^{2}b(d^{2}b+c)}\leq 4d^{2}b+3c.$
\end{clm}
We then find
\begin{equation}
 \int_{M}\chi^{2}|\nabla \phi|^{2}\leq 64n^{4n+2}(1+K)^{2n}\int_{M} \phi^{2} |\nabla \chi|^{2}+3|\phi|_{C^{0}(M)}\cdot |e^{f}-1|_{L^{1}(M)}
\end{equation}


Still let $\rho=\rho_{i}$ and $i \rightarrow \infty$. Same argument as \eqref{equ stripes} and \eqref{equ energy <=0} yields
\begin{equation}\label{equ explicit energy bound}
 \int_{M}|\nabla \phi|^{2}\leq 64n^{4n+2}(1+K)^{2n} c_{vol} |\phi|_{C^{0}(M)}+3|\phi|_{C^{0}(M)}\cdot |e^{f}-1|_{L^{1}(M)}.
\end{equation}
We enlarge the constant to \eqref{equ cD}.
\subsection*{For Theorem \ref{thm}.3}
 By monotone convergence theorem, the established Dirichlet energy bound implies 
$$\lim_{\rho \rightarrow \infty}\int_{\{r\geq \rho\}} |\nabla \phi|^{2}= 0. $$
We apply the argument last section to $\phi-\mu_{i}$, which is also a solution.  Multiplying   $i\partial \overline{\partial}\phi\wedge Q = 0$ by $\chi^{2}(\phi-\mu_{i})$, we still have equality. 
The assumed Poincar\'e inequality \eqref{equ Poincare inequ} and derivation of \eqref{equ stripes} shows 
\begin{eqnarray}\label{equ energy bd using Poincare}
& & \int_{B(\rho_{i})}|\nabla \phi|^{2}=\int_{B(\rho_{i})}|\nabla (\phi-\mu_{i})|^{2}
\leq    \frac{16n^{4n+2}(1+K)^{2n}}{\rho_{i}^{2}}\int_{B(2\rho_{i})\setminus B(\rho_{i})} |\phi-\mu_{i}|^{2} \nonumber
 \\&\leq &  16n^{4n+2}(1+K)^{2n}c_{P}\int_{\{r\geq \underline{\rho}_{i}\}} |\nabla \phi|^{2}\ \textrm{which}\  \longrightarrow 0\ \textrm{as}\ i\rightarrow \infty. 
 \end{eqnarray}
 This again says $\int_{M}|\nabla \phi|^{2}=0$ and $\phi$ is a constant. 

\subsection*{For Corollary \ref{Cor}}
It is obvious from the general uniqueness  \ref{thm} and the volume growth condition. Because the density data $f$  in \cite[Theorem 1.1]{TianYau} and \cite[Proposition 3.4]{HeinThesis} are bounded on the whole $M$ i.e. there exists a positive (finite) number $c_{0}$ such that $|f|\leq c_{0}$,  any  bounded $C^{1,1}-$solution $v$ yields  $\omega_{0}+i\partial\overline{\partial}v$ quasi-isometric to the reference metric $\omega_{0}$ i.e. there is a $c_{1}>0$ possibly depending on $v$ and $f$ such that 
\begin{equation*}
\frac{\omega_{0}}{c_{1}}\leq \omega_{0}+i\partial\overline{\partial}v \leq c_{1}\omega_{0}. 
\end{equation*}
In conjunction with the remark above Corollary \ref{Cor}, if there are two solutions $\phi_{1}$ and $\phi_{2}$, use   $\omega_{\phi_{1}}$ as the new reference metric and denote it by $\omega$,  and denote $\phi_{2}-\phi_{1}$ by $\phi$.  The volume form of $\omega$ and $\omega_{\phi}$ coincide. We apply Theorem \ref{thm}.

 When $\alpha<2$ in the $(K,\alpha,\beta)-$condition \cite[Definition 1.1]{TianYau}, or $\beta<2$ in Hein's $SOB(\beta)-$condition, by quasi-isometry, $\omega$ has  strict sub-quadratic volume growth. Then Theorem \ref{thm}.1 yields the result.

We  elaborate more for  $SOB(2)$.  We do not know whether the solution $\omega$ is $SOB(2)$, though $\omega_{0}$ is by assumption.  Nevertheless,    
 \textit{the  interior Neumann Poincare inequality \cite[Proposition 3.4]{HeinThesis} for $SOB(2)-$reference metric $\omega_{0}$ and the   quasi-isometry}
 $$\frac{\omega}{c_{2}}\leq \omega_{\phi} \leq c_{2}\omega$$
 \textit{still implies the weak Neumann Poincare inequality} \eqref{equ Poincare inequ} \textit{for $\omega$}. 
  Namely, fix a single base point $o$ for both $\omega_{0}$ and $\omega$.  For large enough $c_{3}$ independent of $i$ such that the ball $B(4\rho_{i})$ with respect to $\omega$ is contained in the ball $B(c_{3}\rho_{i})$ with respect to $\omega_{0}$, and the ball $B(\frac{\rho_{i}}{2})$ with respect to $\omega$ contains  the ball $B(\frac{\rho_{i}}{c_{3}})$ with respect to $\omega_{0}$,  we assign the  data 
 $$r_{1}=\frac{\rho_{i}}{c_{3}},\ s=\frac{\rho_{i}}{1000c_{3}},\ r_{2}=c_{3}\rho_{i},\ \kappa=0$$
on  radius and other to  the ball and annuli  in  \cite[Proposition 3.4]{HeinThesis} for the $SOB(2)-$reference metric $\omega_{0}$. 
When $i$ is large,  \cite[(3.4)]{HeinThesis}  implies
\begin{eqnarray}\label{eqn SOB2}
& &\int_{B(2\rho_{i})\setminus B(\rho_{i})}|h-h_{A(\frac{\rho_{i}}{c_{3}}, c_{3}\rho_{i})}|^{2}\omega^{n}\nonumber
\\&\leq & c^{\prime}_{3}\int_{B(c_{3}\rho_{i})\setminus B(\frac{\rho_{i}}{c_{3}})}|h-h_{A(\frac{\rho_{i}}{c_{3}}, c_{3}\rho_{i})}|^{2}\omega^{n}_{0}\leq c_{4}\rho^{2}_{i}\int_{\{r\geq \frac{\rho_{i}}{10c_{3}} \}}|\nabla_{\omega_{0}} f|^{2}\omega^{n}_{0}\nonumber
\\&\leq &c_{5}\rho^{2}_{i}\int_{\{r\geq \frac{\rho_{i}}{c_{5}} \}}|\nabla f|^{2}\omega^{n},
\end{eqnarray}
where the large enough positive constants  $c^{\prime}_{3},\ c_{4},\ c_{5}$ are independent of $\rho_{i}$ or the arbitrary twice differentiable function $h$, and $h_{A(\frac{\rho_{i}}{c_{3}}, c_{3}\rho_{i})}$ is the $\omega^{n}_{0}-$average of $h$ on the closed annulus therein. Note $h_{A(\frac{\rho_{i}}{c_{3}}, c_{3}\rho_{i})}$ is also applied in the first line
of \eqref{eqn SOB2} in the integration against  $\omega^{n}$.  Related to $c_{4}$, the volume ratio $N$ in \cite[(3.4)]{HeinThesis}  is bounded (from above) by the volume constants in $SOB(2)$-condition \cite[Definition 3.1]{HeinThesis} (including the Ricci lower bound).
  This fulfills requirement \eqref{equ Poincare inequ} because $\frac{\rho_{i}}{c_{5}}$ still approaches   $\infty$ as $i\rightarrow \infty$. Then apply Theorem \ref{thm}.3. 
  
  .

\end{document}